\documentclass[graybox]{svmult}

\usepackage{mathptmx}       % selects Times Roman as basic font
\usepackage{helvet}         % selects Helvetica as sans-serif font
\usepackage{courier}        % selects Courier as typewriter font
\usepackage{type1cm}        % activate if the above 3 fonts are
\usepackage{makeidx}         % allows index generation
\usepackage{graphicx}        % standard LaTeX graphics tool
\usepackage{multicol}        % used for the two-column index
\usepackage[bottom]{footmisc}% places footnotes at page bottom

\usepackage{amsfonts}

\makeindex             % used for the subject index
\def\mathsf{\bf}

\def\d{\mathrm d}

\def\E{\mathrm E}
\def\P{\mathrm P}

\def\text{\mbox}

\def\text{\mbox}

\begin{document}
\title*{Modification of Moment-Based Tail Index Estimator: Sums versus Maxima
}
% Use
\titlerunning{Modification of Moment-Based Tail Index Estimator}
%for an abbreviated version of
\author{Natalia Markovich and Marijus Vai\v{c}iulis}
\institute{Natalia Markovich \at V.A. Trapeznikov Institute of Control Sciences of Russian Academy of Sciences,
        Moscow, 117997 Russia, \email{nat.markovich@gmail.com}
\and Marijus Vai{\v c}iulis \at Institute of Mathematics and Informatics, Vilnius University, Akademijos st. 4, LT-08663 Vilnius, Lithuania \email{marijus.vaiciulis@gmail.com}
}
\maketitle
\abstract*{
% 10-15  lines
In this paper we continue the investigation of the SRCEN estimator of the extreme value index $\gamma$ (or the tail index $\alpha=1/\gamma$) proposed in \cite{MCE} for $\gamma>1/2$. We
propose  a new estimator based on the local maximum. This, in fact, is a modification of the SRCEN
estimator to the case $\gamma>0$.
We establish the consistency and asymptotic normality of the newly proposed estimator for i.i.d. data. Also,
a short discussion on the comparison of the estimators is included.
}

\abstract{In this paper, we continue the investigation of the SRCEN estimator of the extreme value index $\gamma$ (or the tail index $\alpha=1/\gamma$) proposed in \cite{MCE} for $\gamma>1/2$. We
propose  a new estimator based on the local maximum. This, in fact, is a modification of the SRCEN
estimator to the case $\gamma>0$.
We establish the consistency and asymptotic normality of the newly proposed estimator for i.i.d. data. Additionally,
a short discussion on the comparison of the estimators is included.
}
%MSC: 62F12, 62G32, 60F05.

\section{Introduction and main results}
Let $X_k, \ k \ge 1$ be non-negative independent, identically
distributed (i.i.d.) random variables (r.v.s) with the
 distribution function (d.f.) $F$.
Suppose that   $F$ belongs to  the domain of attraction of the
Fr\'{e}chet distribution
$$
    \Phi_{\gamma}(x)=\left\{
                       \begin{array}{ll}
                         0, & x\le 0, \\
                         \exp\{-x^{-1/\gamma}\},& x>0,
                       \end{array}
                     \right.
\quad \Phi:=\Phi_1,
$$
which means that there exists normalizing constants $a_m>0$ such
that
\begin{equation}\label{aa02}
    \lim_{m \to \infty} \P\left(\frac{L_{m}}{a_m} \le x \right)=
    \lim_{m \to \infty} F^m\left(a_m x \right)=\Phi_{\gamma}(x),
\end{equation}
for all $x>0$, where $L_{u,v}=\max\{X_{u}, \dots, X_{v}\}$ for $1  \le u \le v$ and  $L_v=L_{1,v}$. The parameter $\gamma>0$ is referred to as positive extreme-value index in the statistical
literature.

Meerschaert and Scheffler \cite{Meer} introduced  the estimator for
$\gamma\ge 1/2$, which is based on the growth rate of the logged
sample variance  of $N$ observations $X_1, \dots, X_N$:
$$
\hat{\gamma}_N=\frac{1}{2\ln(N)} \ln_{+} \left( N s_N^2\right),
$$
where $s_N^2=N^{-1} \sum_{i=1}^N \left( X_i- \bar{X}_N\right)^2$,  $\bar{X}_N=(X_1+\dots+X_N)/N$ and $\ln_+(x)= 0 \vee \ln x$.
%In \cite{Meer}
%is proved that for data in the domain of attraction of
%a stable law the  estimator $\hat{\gamma}_N$ is asymptotically log stable with $\ln(N)$ rate of convergence.

McElroy and Politis \cite{MCE} divided the observations $X_1, \dots, X_N$ into non-intersecting blocks
$\{X_{(k-1)m^2+1}, \dots, X_{km^2}\}$, $1\le k \le [N/m^2]$ of the width $m^2$, while each such block was divided into non-intersecting sub-blocks
of the width $m$. %For the estimation of the parameter
To estimate $\gamma>1/2$ %in the $k$-th block %McElroy and Politis \cite{MCE} considered
%the "centered" (and based on the rate of divergence
%of the logged second moment)  estimator
the  so-called SRCEN  estimator %for  $\gamma>1/2$
was proposed as the sample mean over all blocks:
$$
\hat{\gamma}_N^{(1)}(m)= \frac{1}{[N/m^2]} \sum_{i=1}^{[N/m^2]}\xi_i(m), %\hat{\gamma}_{N,k}^{(1)}(m).
$$
where
\begin{equation}\label{iv000}
%\hat{\gamma}_{N,k}^{(1)}(m)
\xi_i(m) = \frac{\ln\left(
\sum_{j=(i-1)m^2+1}^{im^2} X_j^2\right)}{2\ln(m)}-\frac{1}{m}\sum_{k=1}^m
\frac{\ln\left( \sum_{j=(k-1)m^2+(k-1)m+1}^{(k-1)m^2+k m}
X_j^2\right)}{2\ln(m)},
\end{equation}
and $[\cdot]$ denotes the integer part.
In %the  practical
applications a simple heuristic rule for the choice of sub-block width $m=[N^{1/3}]$, provided in \cite{MCE}, works quite well, see the Monte-Carlo simulation
studies in \cite{MCE}, \cite{Vaic0} and \cite{Vaic}.

Using the inequality of arithmetic and geometric means we obtain that for sample %any observations
$X_1, \dots, X_N$,  $\hat{\gamma}_N^{(1)}(m) \ge 1/2$ holds
with equality if and only if $X^2_{(i-1)m^2+1}=\dots=X^2_{im^2}$, $1 \le i \le [N/m^2]$.
% as far as $\hat{\gamma}_N^{(2)}(m) > 0$ holds.

In this paper we provide an estimator similar to the  SRCEN estimator but one that can
be used for $\gamma>0$, not only for $\gamma>1/2$. Namely, we replace the sums %over $k$
in (\ref{iv000}) by corresponding maxima and introduce the new estimator
$$
\hat{\gamma}_N^{(2)}(m)= \frac{1}{ [N/m^2]} \sum_{i=1}^{[N/m^2]}\widetilde{\xi_i}(m) %\hat{\gamma}_{N,k}^{(2)}(m),
$$
where
$$
%\hat{\gamma}_{N,k}^{(2)}(m)
\widetilde{\xi_i}(m)
= \frac{\ln\left( L_{(i-1)m^2+1, im^2}\right)}{\ln(m)}-\frac{1}{m}\sum_{j=1}^m \frac{\ln\left( L_{(i-1)m^2+(j-1)m+1,(i-1)m^2+j m}\right)}{\ln(m)}.
$$
In fact, the estimator $\hat{\gamma}_N^{(2)}(m)$ is based on the convergence $\E \ln \left(L_m\right)/\ln(m) \to \gamma$ as $m \to \infty$, which implies
\begin{equation}\label{pagrind}
2 \E \left(\frac{ \ln \left(L_{m^2}\right)}{\ln(m^2)}\right)- \E \left(\frac{ \ln \left(L_{m}\right)}{\ln(m)}\right) \to \gamma, \quad m \to \infty.
\end{equation}
Thus, the estimator $\hat{\gamma}_N^{(2)}(m)$
is nothing else, but a moment-type estimator for the left hand side in (\ref{pagrind}).

Note that  $\hat{\gamma}_N^{(2)}(m)$ as well as
 $\hat{\gamma}_N^{(1)}(m)$ are scale-free, i.e., they do not change when $X_j$ is
replaced by $c X_j$ with $c > 0$.

Typically, the estimators, whose constructions are based on the grouping of the observations into the blocks, are well suited for recursive on-line calculations.
In particular, if $\hat{\gamma}_N^{(1)}(m)=\hat{\gamma}_N^{(1)}(m; X_1, \dots, X_{N})$
 denotes the estimate of  $\gamma$
obtained from observations $X_1, \dots, X_{N}$ and we get the next group of updates $X_{N+1}, \dots, X_{N+m^2}$,
then  we obtain
\begin{eqnarray*}&&
\hat{\gamma}_N^{(1)}(m; X_1, \dots,X_{N+m^2})=\frac{1}{\tilde{N}+1}\sum_{i=1}^{\tilde{N}+1}\xi_i(m)
= \frac{1}{\tilde{N}+1}\left(\tilde{N}\hat{\gamma}^{(1)}_N(m)+ \xi_{\tilde{N}+1}(m)\right),
\end{eqnarray*}
denoting $\tilde{N}=[N/m^2]$.
After getting $L$ additional groups $\{X_{N+(k-1)m^2+1}, \dots, X_{N+km^2}\}$, $k=1,...,L$, we have
\begin{eqnarray*}&&
\hat{\gamma}_N^{(1)}(m; X_1, \dots,X_{N+Lm^2})=\frac{1}{\tilde{N}+L}\sum_{i=1}^{\tilde{N}+L}\xi_i(m)
%\\
%&=&\frac{\ell}{(\ell+L)\ell}\left(\sum_{i=1}^{\ell}\xi_i(m)+\sum_{i=\ell+1}^{\ell+L}\xi_{\ell+1}(m)\right)
\\
&=& \frac{1}{\tilde{N}+L}\left(\tilde{N} \hat{\gamma}^{(1)}_N(m)+ \xi_{\tilde{N}+1}(m)+...+\xi_{\tilde{N}+L}(m)\right).
\end{eqnarray*}
It is important  that $\hat{\gamma}_N^{(1)}(m; X_1, \dots,X_{N+Lm^2})$ is obtained using $\hat{\gamma}_N^{(1)}(m)$ after $O(1)$ calculations.
The same is valid for $\widehat{\gamma}^{(2)}_N(m)$ substituting $\xi_i(m)$ by $\widetilde{\xi}_i(m)$.
The discussion on on-line
estimation of the parameter $\gamma>0$ can be found in Section 1.2.3 of \cite{Mark}.

There are situations  when data can be divided
naturally into blocks but only the largest observations within blocks (the block-maxima) are available.
Several such examples are mentioned in \cite{Qi}, see also \cite{Cede}, where battle deaths in
major power wars  between 1495 and 1975 were analyzed.
Then the estimator $\hat{\gamma}_N^{(2)}(m)$ can be applied while the estimators $\hat{\gamma}_N$ and $\hat{\gamma}_N^{(1)}(m)$ are not applicable.

We will formulate our assumptions in terms of a
so-called  quantile function $V$ of the d.f. $F$, which  is defined as the  left continuous generalized inverse:
$$
V(t): =\inf\left\{x \ge 0: \ -\frac{1}{\ln F(x)} \ge t\right\}.
$$
The domain of attraction condition (\ref{aa02}) can be stated in the
following way in terms of $V$:  regarding the d.f. $F$,
(\ref{aa02}) holds if and only if for all $x>0$,
\begin{equation}\label{aa04}
    \lim_{t\to \infty} \frac{V(tx)}{V(t)}=x^{\gamma},
\end{equation}
i.e. the function $V$ varies regularly at infinity with the index  $\gamma>0$ (written $V \in RV_{\gamma}$), see, e.g.,  \cite[p.34]{Haanbook}.

First our result states that $\hat{\gamma}_N^{(2)}(m)$ is a weakly
consistent estimator for  $\gamma>0$. For the sake of completeness we include a corresponding result (as a direct consequence of Prop. 1 in \cite{MCE}) for
the SRCEN estimator $\hat{\gamma}_N^{(1)}(m)$.
\begin{theorem} \label{thm1}
Let observations $X_1,\dots,X_N$ be i.i.d. r.v.s with d.f. $F$.

(i) Suppose $F$ satisfies the first-order condition (\ref{aa04})  with $\gamma>1/2$. Suppose, in addition,
that the probability density function $p(x)$ of $F$ exists and is
bounded, and also that $p(x)/x$ is bounded in a neighborhood of zero. Then for the sequence $m=m(N)$  satisfying
\begin{equation}\label{iv01}
m(N) \to \infty, \quad \frac{N \ln^2 m}{m^2}\to \infty, \quad N\to \infty,
\end{equation}
it holds
\begin{equation}\label{iv01-1e}
\hat{\gamma}_N^{(1)}(m) \ {\buildrel \rm P \over \to} \ \gamma,
\end{equation}
where ${\buildrel \rm P \over \to}$ denotes convergence in probability.

(ii)  Suppose $F$ satisfies (\ref{aa04})  with $\gamma>0$. Suppose, in addition,
\begin{equation} \label{thm1-a}
F(\delta)=0
\end{equation}
for some $\delta>0$.
Then for the sequence $m=m(N)$  satisfying (\ref{iv01}) it holds
\begin{equation}\label{iv01-1}
\hat{\gamma}_N^{(2)}(m) \ {\buildrel \rm P \over \to} \ \gamma.
\end{equation}
\end{theorem}

As usual, in order to get asymptotic normality for estimators the so-called second-order regular variation
condition in some form is assumed.
We recall that the function $V$  is
said to satisfy the second-order condition if for some
 measurable
function $A(t)$ with the constant sign near infinity, which is not identically zero, and $A(t) \to 0$
as $t \to \infty$,
\begin{equation}\label{iv01-1a}
    \lim_{t \to \infty} \frac{\frac{V(tx)}{V(t)}-x^{\gamma}}{A(t)}=x^{\gamma} \frac{x^{\rho}-1}{\rho}
\end{equation}
holds for all $x > 0$ with $\rho < 0$, %$\rho \le 0$,
 which is  a second order parameter.
The function $A(t)$ measures the rate of convergence of $V(tx)/V(t)$
towards $x^{\gamma}$ in (\ref{aa04}), and $|A(t)| \in RV_{\rho}$, see  \cite{Geluk}.

In this paper, we  assume a second order condition stronger than (\ref{iv01-1a}). Namely, we assume that we are in Hall's class of models (see \cite{Hall1}), where
\begin{equation}\label{iv01-1b}
    V(t)=C t^{\gamma}\left(1+\rho^{-1} A(t)\left(1+o(1)\right)\right), \quad t \to \infty
\end{equation}
with $A(t)=\gamma \beta t^{\rho}$, where $C>0$, $\beta \in \mathbb{R}
\setminus \{0\}$ and $\rho<0$. The relation (\ref{iv01-1b}) is
equivalent to
\begin{equation}\label{iv01-1b1}
    F(x)=\exp\left\{ -\left( \frac{x}{C}\right)^{-1/\gamma} \left( 1+\frac{\beta}{\rho} \left( \frac{x}{C}\right)^{\rho/\gamma} +o\left( x^{\rho/\gamma}\right) \right)\right\},
 \quad x \to \infty.
\end{equation}

\begin{theorem} \label{thm2}
Let the observations $X_1,\dots,X_N$ be i.i.d. r.v.s with d.f. $F$.

(i) Suppose $F$ satisfies the second-order condition (\ref{iv01-1b1}) with $\gamma>1/2$ and,
in addition,
that the probability density function $p(x)$ of $F$ exists and it is
bounded, and also that $p(x)/x$ is bounded in a neighborhood of zero. Then for the sequence $m=m(N)$  satisfying
$m \to \infty$ and
\begin{eqnarray*}
  N^{1/2} m^{-2 \vee (-1+\rho) \vee (-2\gamma)} \ln(m) \to 0, \quad {\rm if} \ -1\vee \rho \not = 1-2\gamma, \\
N^{1/2} m^{-2\gamma} \ln^2 (m) \to 0, \quad {\rm if} \ -1\vee \rho  = 1-2\gamma,
\end{eqnarray*}
\begin{equation}\label{norm1}
\frac{N^{1/2} \ln(m)}{m}\left( \hat{\gamma}_N^{(1)}(m) - \gamma\right) \ {\buildrel \rm d \over \to} \ \mathcal{N}\left(0, \frac{\left(\gamma^2-(1/4)\right) \pi^2}{6} \right), \quad N \to \infty,
\end{equation}
holds, where ${\buildrel \rm d \over \to}$ stands for the convergence in distribution.

(ii) Suppose $F$ satisfies (\ref{thm1-a}) and  (\ref{iv01-1b1}) with $\gamma>0$.
Then, for the sequence $m=m(N)$ satisfying (\ref{iv01}) and
\begin{equation}\label{iv01-1c}
    \frac{N^{1/2}}{m} A(m) \to \nu \in (-\infty, +\infty),
\end{equation}
it follows
\begin{equation}\label{ivv01}
\frac{N^{1/2} \ln(m)}{m}\left( \hat{\gamma}_N^{(2)}(m) - \gamma\right) \ {\buildrel \rm d \over \to} \mathcal{N}\left(-\frac{\nu \Gamma(1-\rho)}{\rho}, \frac{\gamma^2 \pi^2}{6} \right), \quad N \to \infty.
\end{equation}
\end{theorem}
The rest of the paper is organized as follows. In the next section we investigate the asymptotic
mean squared error (AMSE) of the introduced estimator, and compare this estimator with several classical estimators, using the same methodology as in \cite{Peng}.
The last section contains the proofs of the results.

\section{Comparison}

The AMSE of the estimator $\hat{\gamma}_N^{(2)}(m)$ is given by
\begin{equation}\label{ivv2}
{\rm AMSE}\left( \hat{\gamma}_N^{(2)}(m) \right):= \frac{1}{\ln^2(m)}\left\{ \frac{\Gamma^2(1-\rho) A^2(m)}{\rho^2}+\frac{\gamma^2\pi^2 m^2}{6N}\right\}.
\end{equation}
Regular variation theory, provided in \cite{Dek} (see also \cite{Peng}), allows us to perform the minimization of the sum in the curly brackets of
(\ref{ivv2}). Namely, under the choice
$$
\bar{m}(N)=\left(\frac{6 \Gamma^2(1-\rho)\beta^2}{-\rho \pi^2}\right)^{1/(2(1-\rho))} N^{1/(2(1-\rho))}\left(1+o(1)\right), \quad N \to \infty,
$$
we have
$$
{\rm AMSE}\left( \hat{\gamma}_N^{(2)}(\bar{m}) \right) \sim \Gamma^2(-\rho)\beta^2\left( \frac{6\beta^2 \Gamma^2(1-\rho)}{\pi^2(-\rho)}\right)^{1/(1-\rho)} \frac{N^{\rho/(1-\rho)}}{\ln^2(N)},\quad N \to \infty.
$$

Probably, the Hill's estimator
$$
\gamma^{(H)}_{N}(k)=\frac{1}{k} \sum_{j=0}^{k-1} \ln\left(\frac{X_{N-j,N}}{X_{N-k,N}} \right),
$$
is the most popular, \cite{Hill}. Here,  $1 \le k \le N$ is a tail
sample fraction, while $X_{1,N} \le X_{2,N} \le \dots \le X_{N,N}$
are order statistics from a sample $X_1, \dots, X_N$. Let us denote
$r=-1 \vee \rho$ and
$$
\upsilon=
\left\{
  \begin{array}{ll}
    \beta, & \hbox{$-1<\rho<0$, } \\
    \beta+(1/2), & \hbox{$\rho=-1$,} \\
    1/2, & \hbox{$\rho<-1$.}
  \end{array}
\right.
$$
From \cite{Peng} it follows that the minimal AMSE of the Hill's estimator under assumption (\ref{iv01-1b1}) satisfies the relation
$$
{\rm AMSE}\left( \gamma_{N}^{(H)}\left( \bar{k} \right) \right) \sim \frac{1-2r}{-2r}\left( \frac{-2r \upsilon^2 \gamma^{2-4r} }{(1-r)^2}\right)^{1/(1-2r)} N^{2r/(1-2r)}, \quad N \to \infty,
$$
where
$$
\bar{k}(N)=\left( \frac{(1-r)^2}{-2r \upsilon^2}\right)^{1/(1-2r)}
N^{-2r/(1-2r)}\left( 1+o(1)\right), \quad N \to \infty.
$$
Now we can compare the estimators $\hat{\gamma}_N^{(2)}(\tilde{m})$ and $\gamma^{(H)}_{N}\left( \bar{k}\right)$.
Denote the relative minimal AMSE in the same way as in \cite{Peng}:
$$
{\rm RMAMSE}(\gamma, \beta, \rho)=\lim_{N\to \infty} \frac{{\rm AMSE}\left( \gamma_{N}^{(H)} \left( \bar{k} \right)\right)}{{\rm AMSE}\left( \hat{\gamma}_N^{(2)}(\bar{m}) \right) }.
$$
Following \cite{Peng} we may conclude that  $\gamma_{N}^{(H)}\left( \bar{k} \right)$ dominates  $\hat{\gamma}_N^{(2)}(\bar{m})$ at the point $(\gamma, \beta, \rho)$ if
${\rm RMAMSE}(\gamma, \beta, \rho)<1$ holds. Note that ${\rm RMAMSE}(\gamma, \beta, \rho)=0$ holds for $-2<\rho<0$, i.e.
 $\gamma_{N}^{(H)}\left( \bar{k} \right)$ dominates  $\hat{\gamma}_N^{(2)}(\bar{m})$, while
for $\rho \le -2$ we have ${\rm RMAMSE}(\gamma, \beta, \rho)=\infty$
and thus,  $\hat{\gamma}_N^{(2)}(\bar{m})$ outperforms
$\gamma_{N}^{(H)}\left( \bar{k} \right)$ in this region of the
parameter $\rho$. It is worth to note that the same conclusion holds
if we replace Hill's estimator by another estimator investigated in
\cite{Peng}.

Unfortunately, it is impossible to compare the
performance of  $\hat{\gamma}_N^{(1)}(m)$ and other estimators taking the AMSE as a measure. By taking $\nu=0$ in (\ref{ivv01})
one can compare the estimators $\hat{\gamma}_N^{(1)}(m)$ and $\hat{\gamma}_N^{(2)}(m)$ under the same block width $m^2$. By comparing variances in the
limit laws  (\ref{norm1}) and (\ref{ivv01}) we conclude that  $\hat{\gamma}_N^{(1)}(m)$ outperforms  $\hat{\gamma}_N^{(2)}(m)$
for $\gamma >1/2$.

\section{Proofs}

Let us firstly provide  preliminary results that are useful in our proofs.

\begin{lemma} \label{lem1}
Let  $X_1,\dots,X_N$ be i.i.d. r.v.s with d.f. $F$. Suppose $F$ satisfies (\ref{aa04})  with $\gamma>0$ and (\ref{thm1-a}). Then
\begin{eqnarray}
\label{prel01}
&&\lim_{m \to \infty} \E \ln \left(\frac{L_{m}}{V(m)} \right) = \chi \gamma, \\
\label{prel02}
&& \lim_{m \to \infty} \E \ln^2 \left(\frac{L_{m}}{V(m)}\right) = \gamma^2 \left( \chi^2+\frac{\pi^2}{6}\right), \\
\label{prel03}
&&\lim_{m \to \infty} \E \ln^4 \left(\frac{L_{m}}{V(m)}\right)  = \gamma^4 \left( \chi^4+\chi^2 \pi^2+\frac{3 \pi^4}{20}+8\chi \zeta(3)\right), \\
\label{prel04}
&&\lim_{m \to \infty} \E \left(\ln \left(\frac{L_{m^2}}{V(m^2)}\right) \ln \left(\frac{L_{m}}{V(m)} \right)\right) = \chi^2 \gamma^2,
\end{eqnarray}
holds, where $\chi\approx 0.5772$ is the Euler–-Mascheroni constant defined
by \\
$\chi=-\int_0^{\infty} \ln(t) \exp\{-t\} \d t$, while $\zeta(t)$ denotes the Riemann zeta function, $\zeta(3) \approx 1.202$.
\end{lemma}
{\it Proof of Lemma \ref{lem1}.} We shall prove (\ref{prel01}).
Let $Y$ be a r.v. with d.f. $\Phi$. It is easy to check that it holds
$$
\ln \left( \frac{L_{m}}{V(m)}\right) \ {\buildrel \rm d \over =} \
 \ln \left( \frac{V(mY)}{V(m)}\right).
$$
%where ${\buildrel \rm d \over =}$ denotes equality in distribution.

By Theorem B.1.9 in \cite{Haanbook}, the assumption $V \in RV_{\gamma}$, $\gamma > 0$  implies  that  for arbitrary
$\epsilon_1>0$, $\epsilon_2>0$ there exists $m_0=m_0(\epsilon_1, \epsilon_2)$ such that for $m \ge m_0$, $my \ge m_0$,
$$
    (1-\epsilon_1) y^\gamma \min\left\{y^{\epsilon_2},y^{-\epsilon_2}\right\} < \frac{V(my)}{V(m)}<(1+\epsilon_1) y^\gamma \max\left\{y^{\epsilon_2},y^{-\epsilon_2}\right\}
$$ holds.
Whence we get that under restriction $0<\epsilon_1<1$ it follows
\begin{equation}\label{Thm01-02}
    \ln (1-\epsilon_1)+(\gamma-  u(y)) \ln(y)  < \ln\left(\frac{V(my)}{V(m)}\right)<\ln(1+\epsilon_1)+ (\gamma+ u(y)) \ln(y),
\end{equation}
where $u(y) = -\epsilon_2 I\{y<1\}+\epsilon_2I\{y \ge1\}$ and $I\{\cdot\}$ denotes the indicator function.

We write for $m>m_0$,
\begin{eqnarray*}
\E\left(\ln \left( \frac{V(m Y)}{V(m)}\right)\right) &=& J_{1,m}+J_{2,m},
\end{eqnarray*}
where
\begin{eqnarray*}
J_{1,m} =\int_0^{m_0/m} \ln\left( \frac{V(my)}{V(m)}\right) \d \Phi(y), \quad
J_{2,m} =\int_{m_0/m}^{\infty} \ln\left( \frac{V(my)}{V(m)}\right) \d \Phi(y).
\end{eqnarray*}
The statement (\ref{prel01}) follows from
\begin{eqnarray} \label{Thm01-03}
\lim_{m \to \infty} J_{1,m} &=&0, \\
\label{Thm01-04} \lim_{m \to \infty} J_{2,m} &=& \chi \gamma.
\end{eqnarray}
Substituting $my=t$ we get
\begin{eqnarray*}
\left|J_{1,m}\right| &\le &\int_0^{m_0} \left|\ln\left( \frac{V(t)}{V(m)}\right)\right| \d \Phi(t/m) \\
&=& \int_0^{m_0} \left|\ln V(t)\right| \d \Phi(t/m) + \Phi(m_0/m) \left|\ln V(m)\right|.
\end{eqnarray*}
By using $\d \Phi(t/m)=m \Phi\left(t/(m-1)\right) \d \Phi(t)$ we obtain
\begin{eqnarray*}
\left|J_{1,m}\right| &\le & m   \Phi\left(m_0/(m-1)\right) \int_0^{m_0} \left|\ln V(t)\right| \d \Phi(t) + \Phi(m_0/m) \left|\ln V(m)\right|.
\end{eqnarray*}
Assumption (\ref{thm1-a}) ensures $V(0) \ge \delta$, which implies $\int_0^{m_0} \left|\ln V(t)\right| \d \Phi(t)<\infty$.
Since the sequence $V (n)$ is of a polynomial growth and  $\Phi(m_0/m)=\exp\{-m/m_0\}$
tends to zero exponentially fast, then relation  (\ref{Thm01-03})  follows.

To prove (\ref{Thm01-04}) we use inequality (\ref{Thm01-02}). Then we obtain
$$
\left|J_{2,m}-\chi \gamma  \right| \le \max\left\{
-\ln(1-\epsilon_1), \ln(1+\epsilon_1)\right\} + \epsilon_2 \E \left|
\ln(Y)\right|+\gamma  \int_0^{m_0/m} \left|\ln\left( y\right)\right|
\d \Phi(y).
$$
One can check that $\E \left| \ln(Y)\right|=\chi-2 {\rm Ei}(-1)$,
where ${\rm Ei}(x)$, $x \in \mathbb{R} \setminus \{0\}$ denotes the
exponential integral function, ${\rm Ei}(-1)\approx -0.219384$.

Since $\epsilon_1>0$ and $\epsilon_2>0$
may be taken arbitrary small, the proof of relation (\ref{Thm01-04}) will be finished if we show that
$\int_0^{m_0/m} \left|\ln\left( y\right)\right| \d \Phi(y) \to 0$, $m \to \infty$.
Substituting $t=my$ we get
\begin{eqnarray*}
\int_0^{m_0/m}\left| \ln\left( y\right)\right| \d \Phi(y) &=& \int_0^{m_0} \left| \ln(t/m)\right| \d  \Phi(t/m) \\
&=& m \int_0^{m_0} \left| \ln(t/m)\right| \Phi(t/(m-1))\d \Phi(t) \\
&\le& m \Phi(m_0/(m-1)) \left( \ln(m)+\E  \left| \ln(Y)\right|\right)\to 0,
\end{eqnarray*}
as $m \to \infty$. This completes the proof of (\ref{Thm01-04}), and also of  relation (\ref{prel01}).

Proofs of relations (\ref{prel02}) and (\ref{prel03}) are similar and  thus are skipped. It remains to prove (\ref{prel04}).
We note that $L_m$ and $L_{m+1, m^2}$ are independent r.v.s and $L_{m^2}=L_m \vee L_{m+1, m^2}$.
Let $Y_1$ and $Y_2$ are independent  r.v.s with d.f.  $\Phi$. Then it holds
$$
\ln \left(\frac{L_{m^2}}{V(m^2)} \right) \ln \left(\frac{L_{m}}{V(m)} \right)\ {\buildrel \rm d \over =} \ \ln \left(\frac{V(mY_1) \vee V(m(m-1)Y_2) }{V(m^2)} \right) \ln \left(\frac{V(mY_1)}{V(m)}\right),
$$
and consequently,
$$
\E \left(\ln \left(\frac{L_{m^2}}{V(m^2)} \right) \ln \left(\frac{L_{m}}{V(m)} \right)\right)
=\E \left(\ln \left(\frac{V(mY_1) \vee V(m(m-1)Y_2) }{V(m^2)} \right) \ln \left(\frac{V(mY_1)}{V(m)}\right)\right).
$$
Let us recall that $V(t)$, $ t \ge 0$ is a non-decreasing function, see, e.g., Prop. 2.3 in \cite{Embr}. By using this property we obtain
$$
\E \left(\ln \left(\frac{V(mY_1) \vee V(m(m-1)Y_2) }{V(m^2)} \right) \ln \left(\frac{V(mY_1)}{V(m)}\right)\right)=J_{3,m}+J_{4,m}+J_{5,m},
$$
where
\begin{eqnarray*}
  J_{3,m} &=&  \E \left(\ln \left(\frac{V(mY_1)}{V(m^2)} \right) \ln \left(\frac{V(mY_1)}{V(m)}\right) I\{ Y_1>(m-1)Y_2\}\right), \\
  J_{4,m} &=&  \E \left(\ln \left(\frac{ V(m(m-1)Y_2)}{V(m^2)} \right)\right) \E \left( \ln \left(\frac{V(mY_1)}{V(m)}\right) \right), \\
  J_{5,m} &=& \E \left(\ln \left(\frac{ V(m(m-1)Y_2)}{V(m^2)} \right) \ln \left(\frac{V(mY_1)}{V(m)}\right) I\{ Y_1>(m-1)Y_2\}\right).
\end{eqnarray*}
Let us rewrite quantity $J_{4,m}$ as follows:
$$
  J_{4,m} =  \left\{\ln\left( \frac{ V(m(m-1))}{V(m^2)} \right)+\E \ln \left(\frac{ L_{m(m-1)}}{V(m(m-1))} \right)\right\}
\E \ln \left(\frac{ L_{m}}{V(m)} \right).
$$
For any $\epsilon>0$ there exists natural $\tilde{m}_0$ such
that $1/m<\epsilon$ for $m \ge m_0$. Then $V(m^2(1-\epsilon))/ V(m^2) \le V\left(
m^2(1-1/m)\right)/V(m^2) \le 1$. By (\ref{aa04}) we get $V(m^2(1-\epsilon))/V(m^2)
\to (1-\epsilon)^{\gamma}$, $m \to \infty$. Since $\epsilon>0$ can be taken arbitrary small,  the relation $V(m(m-1))/V(m^2)\to 1$,
$m\to \infty$ holds. By using the last relation and (\ref{prel01})
we deduce that $J_{4,m} \to \chi^2 \gamma^2$ holds as $m \to \infty$.

Next, we have
\begin{eqnarray*}
  J_{3,m} &=&  \E \left(\ln^2 \left(\frac{V(mY_1)}{V(m)} \right) I\{ Y_1>(m-1)Y_2\}\right) \\
  &&+ \ln \left(\frac{V(m)}{V(m^2)} \right) \E \left(\ln \left(\frac{V(mY_1)}{V(m)} \right) I\{ Y_1>(m-1)Y_2\}\right).
\end{eqnarray*}
We apply the H\"{o}lder's inequality to get
\begin{eqnarray*}
  \left|J_{3,m}\right| &\le&  \left\{\E \ln^4\left( \frac{L_m}{V(m)}\right)\right\}^{1/2} \left\{\P( Y_1>(m-1)Y_2)\right\}^{1/2} \\
  &&+ \left|\ln \left(\frac{V(m)}{V(m^2)} \right)\right|  \left\{\E \ln^2\left( \frac{L_m}{V(m)}\right)\right\}^{1/2} \left\{\P( Y_1>(m-1)Y_2)\right\}^{1/2}.
\end{eqnarray*}
We find that $\P( Y_1>(m-1)Y_2) =1/m$ holds. Let us recall the well-known property of regularly varying functions:
if $V \in RV_{\gamma}$, then
\begin{equation}\label{lem1-04}
    \lim_{m \to \infty} \frac{\ln V(m)}{\ln(m)}=\gamma,
\end{equation}
see, e.g., Prop. B.1.9 in \cite{Haanbook}. By using (\ref{lem1-04}) we obtain $ \ln\left(V(m^2)/V(m)\right) \sim \gamma \ln (m)$, $m \to \infty$.
Thus, keeping in mind (\ref{prel02}) and (\ref{prel03}) we obtain $\left|J_{3,m}\right|=O\left(m^{-1/2} \ln(m) \right)$, $m \to \infty$.
By a similar argument we obtain $\left|J_{5,m}\right|=O\left(m^{-1/2} \right)$, $m \to \infty$. This finishes the proof of (\ref{prel04}) and Lemma \ref{lem1}.

\bigskip

{\it Proof of Theorem \ref{thm1}.} First we prove (\ref{iv01-1}). Let us rewrite
\begin{equation} \label{Thm02-00}
\hat{\gamma}_N^{(2)}(m)= \gamma + \left\{\E\hat{\gamma}_N^{(2)}(m) -\gamma\right\}  +S_N(m),
\end{equation}
where
\begin{eqnarray}
\E \hat{\gamma}_N^{(2)}(m) -\gamma &=& \left\{\frac{\ln V(m^2)-\ln V(m)}{\ln(m)}- \gamma\right\} \nonumber \\
&&+ \frac{1}{\ln(m)} \left(\E \ln\left( \frac{L_{m^2}}{V(m^2)}\right) -\E \ln\left( \frac{L_{m}}{V(m)}\right)\right)
\label{poslq}
\end{eqnarray}
and
\begin{eqnarray*}
S_{N}(m)&=&\frac{1}{[N/m^2]\ln(m)} \sum_{i=1}^{[N/m^2]} \bigg\{\left\{\ln\left( \frac{L_{(i-1)m^2+1, im^2}}{V(m^2)}\right)-\E \ln\left( \frac{L_{m^2}}{V(m^2)}\right)\right\}\\
&& \quad \quad \quad   -
\frac{1}{m}\sum_{j=1}^m \left\{\ln\left( \frac{L_{(i-1)m^2+(j-1)m+1,(i-1)m^2+j m}}{V(m)}\right)-\E \ln\left( \frac{L_{m}}{V(m)}\right)\right\}\bigg\}.
\end{eqnarray*}
By combining  (\ref{prel01}) and (\ref{lem1-04}) we deduce that $\E \hat{\gamma}_N^{(2)}(m) -\gamma \to 0$, $m \to \infty$.  Thus, it is enough to prove
that $S_{N}(m)  \ {\buildrel \rm P \over \to} \ 0$ as $N \to \infty$. By Chebyshev's inequality, for any $\epsilon > 0$ it holds
$\P\left( \left|S_{N}(m)\right|>\epsilon\right) \le \epsilon^{-2} \E\left(S_{N}(m)\right)^2.$
We have
\begin{eqnarray}
&&\E\left(S_{N}(m)\right)^2 = \frac{1}{[N/m^2]\ln^2(m)} \bigg\{ {\rm Var}\left( \ln\left( \frac{L_{m^2}}{V(m^2)}\right)\right) \nonumber \\
&&-2 {\rm Cov} \left(\ln\left( \frac{L_{m^2}}{V(m^2)} \right), \ln\left( \frac{L_{m}}{V(m)} \right)\right)
+\frac{1}{m} {\rm Var} \left(\ln\left( \frac{L_{m}}{V(m)}\right)\right)\bigg\}.
\label{varq}
\end{eqnarray}
Use (\ref{prel01})-(\ref{prel02}) and (\ref{prel04}) to deduce that  the sum in the curly brackets has a finite limit as $m \to \infty$.
Thus, assumption (\ref{iv01}) ensures $\E\left(S_{N}(m)\right)^2 \to 0$, $m \to \infty$. This finishes the proof of (\ref{iv01-1}).

Consider now (\ref{iv01-1e}), where the restriction $\gamma>1/2$ holds.
Assumption  (\ref{aa04}) is equivalent to $1-F \in RV_{-1/\gamma}$. %RV_{-1/\gamma}. % relations V \in RV_{\gamma} and 1-F RV_{-1/\gamma} are equivalent !
By the Representation Theorem (see, Thm. B.1.6. in \cite{Haanbook}), there exists a function $\ell \in RV_0$, such that
\begin{equation} \label{uodega}
1-F(x^{1/2})=x^{-1/(2\gamma)} \ell\left(x^{1/2}\right), \quad x \to \infty.
\end{equation}
Following the Mijnheer Theorem (see, Thm. 1.8.1 in \cite{SamT}), we  determine the norming function $a(m) \in RV_{2\gamma}$ from
\begin{equation}\label{norm1A}
\lim_{m\to \infty} \frac{m \ell \left( a^{1/2}(m) \right)}{\left(a(m)\right)^{1/(2\gamma)}} = d(\gamma), \quad d(\gamma)=%\frac{\Gamma(1-1/\gamma)}{\cos\left( \pi/(4\gamma)\right)}.
\Gamma(1-1/(2\gamma))\cos\left( \pi/(4\gamma)\right).
\end{equation}
Put $Q(m)=(X_1^2+\dots+X_m^2)/a_m$. Then
$ Q(m) \ {\buildrel \rm d \over \to} Z,
$ as $m \to \infty$,
where $Z$ is totally skewed to the right $1/(2\gamma)$-stable r.v. with characteristic function
\begin{equation}\label{Thm2-02a}
\E\exp\{i\theta Z\}=\exp\left\{ -|\theta|^{1/(2\gamma)} \left( 1-
{\rm i} \ {\rm sgn}(\theta) \tan \left( \frac{\pi}{4
\gamma}\right)\right)\right\}.
% Z \sim S_{-1/(2\gamma)} (1, 1,0), see for definition in Samorodnitsky&Taqqu.
\end{equation}
Similarly to (\ref{Thm02-00}) we use the decomposition
$$
\hat{\gamma}_N^{(1)}(m)= \gamma +\left\{\E \hat{\gamma}_N^{(1)}(m)-\gamma\right\}+\tilde{S}_N(m),
$$
where
\begin{eqnarray*}
%\tilde{B}(m) &=& \frac{\ln a(m^2)-\ln a(m)}{2\ln(m)}+ \frac{1}{2\ln(m)} \left(\E \ln\left( \frac{Q(m^2)}{a(m^2)}\right) -\E \ln\left( \frac{Q(m)}{a(m)}\right)\right), \\
\tilde{S}_{N}(m)&=&\frac{1}{2[N/m^2]\ln(m)} \sum_{i=1}^{[N/m^2]} \left\{\ln\left( \sum_{j=(i-1)m^2+1}^{im^2} \frac{X_j^2}{a(m^2)}\right)-\E \ln Q(m^2)\right\}\\
&& \quad \quad \quad  \ \  -
\frac{1}{m}\sum_{j=1}^m \left\{\ln\left( \sum_{j=(i-1)m^2+(i-1)m+1}^{(i-1)m^2+i m}
\frac{X_j^2}{a(m)}\right)-\E \ln Q(m) \right\}.
\end{eqnarray*}
The bias of the estimator $\hat{\gamma}_N^{(1)}(m)$ is given by
$\E \hat{\gamma}_N^{(1)}(m)-\gamma = \Delta (m^2) -(1/2) \Delta(m)$,
where
$$
\Delta(m)= \frac{\ln a(m)}{\ln m}- 2\gamma +\frac{1}{\ln m} \left\{\E \ln Q(m)- \E \ln Z \right\}.
$$
In Prop. 1-2 of \cite{MCE} it is proved
\begin{eqnarray}\label{Thm2-02}
&&\E \ln Q(m) \to \E\ln Z, \quad \E \ln^2 Q(m) \to \E\ln^2 Z,  \\
\label{Thm2-03}
&&{\rm Cov}\left(\ln Q(m^2), \ln Q(m) \right) \to 0, \quad m \to \infty.
\end{eqnarray}
It is worth to note that the moments $\E\ln Z$ and $\E\ln^2 Z$  can
be found explicitly. Indeed, there is a direct connection between
moments of order $r<1/(2\gamma)$ and log-moments of order $k\in
\mathbb{N}$:
\begin{equation}\label{thm02-04}
    E \ln^k Z = \frac{\d^k}{\d r^k} E Z^r \bigg|_{r=0},
\end{equation}
see \cite{Zol}. Regarding the moments $\E Z^r$, the following relation is proved  in Section 8.3 of \cite{Pae}:
\begin{equation}\label{thm02-03}
\E Z^r=\frac{\Gamma(1-2\gamma r)}{\Gamma(1-r)}\left(1+\tan^2\left(\frac{\pi}{4 \gamma}\right) \right)^{\gamma r}, \quad -1<r<1/(2\gamma).
\end{equation}
By using (\ref{thm02-04}) and (\ref{thm02-03}) we obtain
\begin{eqnarray} \label{first}
  \E \ln Z &=& -\chi+2 \chi  \gamma +\gamma  \ln \left(\tan ^2\left(\frac{\pi }{4 \gamma }\right)+1\right),\\
  \E \ln^2 Z &=&\chi^2-\frac{\pi ^2}{6}+ 4 \chi^2 \gamma^2-4 \chi^2 \gamma+\frac{2 \pi ^2 \gamma ^2}{3}+\gamma ^2 \log ^2\left(\tan ^2\left(\frac{\pi }{4 \gamma }\right)+1\right) \nonumber \\
  && +4 \chi  \gamma ^2 \log \left(\tan ^2\left(\frac{\pi }{4 \gamma }\right)+1\right) -2 \chi  \gamma  \log \left(\tan ^2\left(\frac{\pi }{4 \gamma }\right)+1\right).
\label{second}
\end{eqnarray}
%which gives ${\rm Var} \left( \ln Z\right)=(4\gamma^2-1) \pi^2/6$.

We combine (\ref{lem1-04}) and the first relation in (\ref{Thm2-02}) to deduce that $\Delta(m) \to 0$, $m \to \infty$, which implies
$\E \hat{\gamma}_N^{(1)}(m)-\gamma \to 0$, $m\to \infty$. Thus, relation (\ref{iv01-1e}) will be proved if we show that under assumptions (\ref{iv01}),
$\E \left( \tilde{S}_{N}(m)\right)^2 \to 0$. The last relation can be verified by using (\ref{Thm2-02}) and (\ref{Thm2-03}), and
\begin{equation} \label{idisp2}
  \E \left( \tilde{S}_{N}(m)\right)^2 =\frac{{\rm Var}\left( \ln Q(m^2)\right)-2 {\rm Cov} \left\{\ln Q(m^2), \ln Q(m) \right\}+ m^{-1} {\rm Var}\left( \ln Q(m)\right)}{4 [N/m^2] \ln^2(m)}.
\end{equation}
This completes the proof of Theorem \ref{thm1}.

\bigskip

{\it Proof of  Theorem \ref{thm2}.}
In view of decomposition (\ref{Thm02-00}),
the assertion (\ref{ivv01}) follows from
\begin{eqnarray}
\label{thm02-01s}
\E\left(S_{N}(m)\right)^2 &\sim& \frac{\pi^2 \gamma^2 m^2}{6 N \ln^2(m)}, \\
\label{thm02-01}
  \left\{ \E \left( S_{N}(m)\right)^2\right\}^{-1/2} S_{N}(m) &{\buildrel \rm d \over \to}& \mathcal{N}(0,1), \\
\label{thm02-02}
  \frac{N^{1/2} \ln(m)}{m} \left( \E \hat{\gamma}_N^{(1)}(m) -\gamma \right) &\to& -\frac{\nu \Gamma(1-\rho)}{\rho}, \quad N \to \infty,
\end{eqnarray}
where $\nu$ is the same as in (\ref{iv01-1c}).

Relation (\ref{thm02-01s}) follows from (\ref{varq}) by applying (\ref{prel01})-(\ref{prel02}) and (\ref{prel04}).
To prove (\ref{thm02-01}), by using (\ref{prel01})-(\ref{prel04}) we  check the $4$-th order Lyapunov condition
for i.i.d. random variables forming a triangular array. We skip standard details.

By using (\ref{iv01-1b}) we obtain
$$
\frac{\ln V(m)}{\ln(m)}-\gamma = \frac{1}{\ln(m)} \left\{\ln(C) +\frac{A(m)}{\rho}(1+o(1))\right\}, \quad m \to \infty.
$$
Following  the proof of Lemma 2 in \cite{Vaic} one can obtain
$$
\E \ln\left( \frac{L_{m}}{V(m)}\right)-\chi \gamma = \frac{\Gamma(1-\rho)-1}{\rho} A(m)\left(1+o(1)\right), \quad m \to \infty.
$$
We combine the last two relations, assumption (\ref{iv01-1c}) and decomposition (\ref{poslq}) to verify (\ref{thm02-02}).

Let us discuss the proof of (\ref{norm1}) now. Relations (\ref{Thm2-02}), (\ref{Thm2-03}), (\ref{first})-(\ref{idisp2}) imply
$\E \left(\tilde{S}_{N}(m)\right)^2 \sim  m^2 N^{-1} \ln^{-2} (m) \left( \gamma^2-(1/4)\right) \pi^2 / 6 $, $N \to \infty$. In view of the last relation
it is enough to prove that
\begin{eqnarray}
\label{thm02-05}
 \left\{ {\rm Var} \left(\tilde{S}_{N}(m)\right)\right\}^{-1/2} \tilde{S}_{N}(m)  &{\buildrel \rm d \over \to}& \mathcal{N}\left(0,1\right), \\
\label{thm02-06}
 \E \hat{\gamma}_N^{(1)}(m) -\gamma &=& \left\{
  \begin{array}{ll}
    O\left( m^{-1 \vee \rho \vee (1-2\gamma)} \right), & \hbox{$-1 \vee \rho \not = 1-2\gamma$}, \\
    O\left( m^{1-2\gamma} \ln(m)\right), & \hbox{$-1 \vee \rho  = 1-2\gamma$.}
  \end{array}
\right.
\end{eqnarray}
We skip a standard proof of (\ref{thm02-05}) and focus on the
investigation of the bias $\E \hat{\gamma}_N^{(1)}(m) -\gamma$.
Firstly, we  prove that
\begin{equation}\label{biaso01}
\frac{\ln a(m^2) - \ln a(m)}{2\ln(m)}-\gamma=O\left( \frac{m^{-1 \vee\rho}}{\ln(m)}\right), \quad m \to \infty.
\end{equation}
The relation (\ref{iv01-1b1}) can be written in the form $1-F(x)=x^{-1/\gamma} \ell(x)$, $x\to \infty$, where function $\ell \in RV_0$ has the form
\begin{equation} \label{lkf01}
\ell(x)=    C^{1/\gamma} \left(1+ \tilde{C}(\beta, \rho) \left(x/C\right)^{(-1 \vee\rho)/\gamma}+o\left( x^{(-1 \vee\rho)/\gamma}\right)\right), \quad x \to \infty,
\end{equation}
where
$$
\tilde{C}(\beta, \rho)=
\left\{
  \begin{array}{ll}
    \beta/\rho, & \hbox{$-1<\rho<0$,} \\
    -(2\beta-1)/\rho, & \hbox{$\rho=-1$, $\beta \not=1/2$,} \\
    -1/2, & \hbox{$\rho<-1$.}
  \end{array}
\right.
$$
Now, by using (\ref{norm1A}), one can find that under assumption (\ref{iv01-1b1}) the norming function satisfies the asymptotic relation
$$
a(m)=\left(C^{1/\gamma}/d(\gamma)\right)^{2\gamma} m^{2\gamma}
\left(
1+ 2 \gamma \tilde{C}(\beta, \rho) d^{-(-1\vee \rho)}(\gamma) m^{-1 \vee\rho}+o\left( m^{-1 \vee\rho}\right)
\right)
$$
as $m \to \infty$, while the last relation implies (\ref{biaso01}).

We claim that
\begin{equation}\label{bias02}
    \frac{\E\ln Q(m)-\E\ln  Z}{\ln m}= \left\{
  \begin{array}{ll}
    O\left( m^{-1 \vee \rho \vee (1-2\gamma)} \right), & \hbox{$-1 \vee \rho \not = 1-2\gamma$}, \\
    O\left( m^{1-2\gamma} \ln(m)\right), & \hbox{$-1 \vee \rho  = 1-2\gamma$}
  \end{array}
\right.
\end{equation}
as $m \to \infty$.

Then terms $\ln^{-1}(m^2) \left\{ \E\ln Q(m^2)-\E\ln  Z \right\}$ and $\left( 2\ln(m)\right)^{-1}\left\{\ln a(m^2)-\ln a(m)\right\}-\gamma$ are negligible
with respect to $\ln^{-1}(m) \left\{ \E\ln Q(m)-\E\ln  Z \right\}$ and thus, the relation (\ref{thm02-06}) follows.

To verify (\ref{bias02}) we use the similar
decomposition $\E\ln Q(m)-\E\ln  Z =R_{1,m}-R_{2,m}-R_{3,m}$ as in  the proof of Prop. 3 in \cite{MCE}, where
\begin{eqnarray*}
  R_{1,m} &=&  \int_0^{\infty} \left\{\P\left(\ln Q(m) >x \right) -\P \left(\ln Z>x \right)\right\} \d x, \\
  R_{2,m} &=& \int_{-\ln m}^0 \left\{\P\left(\ln Q(m)<x \right) -\P \left(\ln Z <x \right)\right\} \d x, \\
  R_{3,m} &=& \int_{-\infty}^{-\ln m} \left\{\P\left(\ln Q(m)<x \right) -\P \left(\ln Z <x \right)\right\} \d x.
\end{eqnarray*}
By using substitution  $t=\exp\{x\}$ we obtain
$$
R_{1,m} =  \int_1^{\infty} t^{-1} \left\{\P\left(Q(m)>t \right) -\P \left( Z>t \right)\right\} \d t.
$$
Similarly we get
$R_{2,m} = \int_{1/m}^1 t^{-1}\left\{\P\left(Q(m)<t \right) -\P \left(Z <t \right)\right\} \d t$.
%After an adaptation of Proposition in \cite{Daugavet} to our purposes,
From  Corollary 2 in \cite{Daugavet} it follows
\begin{eqnarray*}\label{bias01}
    &&\sup_{t \ge 0} f_{\gamma}(t) \left| \P\left(Q(m)>t \right) -\P \left( Z>t \right)\right| =O\left( \lambda\left( m^{2\gamma}\right)+m^{-2\gamma}\right), \quad m \to \infty,
\end{eqnarray*}
where $f_{\gamma}(t)= 1 + t^{2\gamma} \ln^{-2}\left( {\rm e} +t\right)$ and $\lambda(R)=\lambda_1(R)+R^{-1+1/(2\gamma)} \lambda_2(R)$, $R>0$, where
\begin{eqnarray*}
  \lambda_1(R) &=& \sup_{u \ge R}  u^{1/(2\gamma)} \left| \P(X_1^2>u)- \P \left( Z> u \right)\right|, \\
  \lambda_2(R) &=& \int_0^R \left| \P(X_1^2>u)- \P \left( Z> u \right)\right| \d u.
\end{eqnarray*}
It is well-known that $P(Z>x)= C_1 x^{-1/(2\gamma)}\left( 1+C_2 x^{-1/(2\gamma)}+o\left(x^{-1/(2\gamma)}\right)\right)$, $x \to \infty$ holds, where $C_k=C_k(\gamma)$ are some
constants. The asymptotic of $\P(X_1^2>u)$ is given in (\ref{uodega}), where   a function $\ell$ slowly varying at infinity is given in (\ref{lkf01}).
Recall that $\hat{\gamma}_N^{(1)}(m)$ is a scale-free estimator. Thus, without loss of generality, we may assume that the scale parameter $C$ in (\ref{lkf01}) satisfies
$C^{1/\gamma}=C_1$. Then we have
\begin{equation} \label{h1}
\P(X_1^2>x)- \P \left( Z> x \right)=D x^{(-2 \vee (\rho-1))/(2\gamma)}+o \left( x^{(-2 \vee (\rho-1))/(2\gamma)}\right), \quad x \to \infty,
\end{equation}
where $D \not =0$ is some constant. By applying (\ref{h1}) we obtain
immediately $\lambda_1\left(m^{2\gamma}\right)=O\left( m^{-1 \vee
\rho}\right)$, $m \to \infty$. If  $-2 \vee (\rho-1)>-2\gamma$,
by  ex. 1.2  in \cite{Fed}, a relation $f(x) \sim x^r$, $x \to \infty$ implies
\begin{equation} \label{Fedasimpt}
\int_0^x f(t) \d t \sim
\left\{
  \begin{array}{ll}
    x^{r+1}/(r+1), & \hbox{$r>-1$,} \\
    \ln(x), & \hbox{$r=-1$,}
  \end{array}
\right.
\quad x \to \infty
\end{equation}
and thus we obtain $m^{1-2\gamma}\lambda_2(m^{2\gamma})=O\left( m^{-1 \vee
\rho}\right)$, $m \to \infty$. In the case  $-2 \vee
(\rho-1)=-2\gamma$, by applying (\ref{Fedasimpt}) one more time
we get $m^{1-2\gamma}\lambda_2(m^{2\gamma})=O\left( m^{1-2\gamma}
\ln(m)\right)$, $m \to \infty$. As for the case $-2 \vee
(\rho-1)<-2\gamma$, we have
$m^{1-2\gamma}\lambda_2(m^{2\gamma})=O\left( m^{1-2\gamma} \right)$,
$m \to \infty$. By putting the obtained results together we get
$$
\sup_{t \ge 0} f_{\gamma}(t) \left| \P\left(Q(m)>t \right) -\P \left( Z>t \right)\right| =
\left\{
  \begin{array}{ll}
    O\left( m^{-1 \vee \rho \vee (1-2\gamma)}\right), & \hbox{$-1 \vee \rho \not = 1-2\gamma$}, \\
    O\left( m^{1-2\gamma} \ln(m)\right), & \hbox{$-1 \vee \rho  = 1-2\gamma$}
  \end{array}
\right.
$$
as $m \to \infty$.

Applying the last asymptotic relation  we obtain immediately
$$
%  \left|R_{1,m}\right| &=& O\left( m^{-1 \vee \rho \vee -2\gamma}\right)+ O\left( m^{1-2\gamma} \ln (m)\right), \\
  \left|R_{2,m}\right|
  =\left\{
  \begin{array}{ll}
    O\left( m^{-1 \vee \rho \vee (1-2\gamma)} \ln(m)\right), & \hbox{$-1 \vee \rho \not = 1-2\gamma$}, \\
    O\left( m^{1-2\gamma} \ln^2(m)\right), & \hbox{$-1 \vee \rho  = 1-2\gamma$}
  \end{array}
\right.
$$
and $\left|R_{1,m}\right|=o\left( \left|R_{2,m}\right|\right)$ as $m \to \infty$. %As for the term $R_{3,m}$,
Since the relation $|R_{3,m}|=O\left(m^{-1}\right)=o\left( \left|R_{2,m}\right|\right)$, $m \to \infty$ holds (see proof of Prop. 3 in \cite{MCE}), %see  \cite{MCE}.
the statement of Theorem \ref{thm2} follows. %is proved.

\end{document}